\documentclass[11pt]{article}
\usepackage{indentfirst}
\usepackage{amssymb}
\usepackage{epsfig}

\def\be{\begin{equation}}
\def\ee{\end{equation}}
\def\bea{\begin{eqnarray}}
\def\eea{\end{eqnarray}}
\def\bes{\begin{eqnarray*}}
\def\ees{\end{eqnarray*}}

\def\nn{\nonumber}
\def\<{\langle}
\def\>{\rangle}
\def\lb{\label}
\def\bs{\setminus}

\def\R{{\bf R}}
\def\C{{\bf C}}
\def\Z{{\bf Z}}
\def\N{{\bf N}}
\def\U{{\bf U}}
\def\Q{{\bf Q}}
\def\T{{\bf T}}

\def\aa{{\alpha}}

\def\ga{{\gamma}}

\def\th{{\theta}}
\def\om{{\omega}}

\def\ep{{\epsilon}}

\def\sg{{\sigma}}

\def\Sg{{\Sigma}}

\def\H{{\cal H}}
\def\T{{\cal T}}

\def\P{{\cal P}}

\def\I{{\cal I}}

\def\Sp{{\rm Sp}}
\def\mod{{\rm mod}}
\def\dm{{\rm \diamond}}

\def\hb{\vrule height0.18cm width0.14cm $\,$}

\title{Multiplicity and stability of closed characteristics on compact convex P-cyclic symmetric hypersurfaces in ${\bf R}^{2n}$}
\author{Hui Liu\thanks{Partially supported by NSFC (Nos. 11771341, 12022111).
E-mail: huiliu00031514@whu.edu.cn}\\
School of Mathematics and Statistics, Wuhan University,
\\Wuhan 430072, Hubei, P. R. China }
\date{}

\begin{document}

\maketitle

\begin{abstract}
{\it  Let $\Sigma$ be a compact convex hypersurface in ${\bf R}^{2n}$ which is P-cyclic symmetric, i.e., $x\in \Sigma$ implies $Px\in\Sigma$
with P being a $2n\times2n$ symplectic orthogonal matrix and satisfying $P^k=I_{2n}$, $ker(P^l-I_{2n})=0$ for $1\leq l< k$, where $n, k\geq2$.
In this paper, we prove that there exist at least $n$ geometrically distinct closed characteristics on $\Sigma$, which solves a longstanding
conjecture about the multiplicity of closed characteristics for a broad class of compact convex hypersurfaces with symmetries(cf.,Page 235 of \cite{Eke1}).
Based on the proof, we further prove that if the number of geometrically distinct closed characteristics on $\Sigma$ is finite,
then at least $2[\frac{n}{2}]$ of them are non-hyperbolic; and if the number of geometrically distinct
closed characteristics on $\Sigma$ is exactly $n$ and $k\geq3$, then all of them are P-cyclic symmetric, where a closed characteristic $(\tau, y)$ on $\Sigma$ is called
P-cyclic symmetric if $y({\bf R})=Py({\bf R})$. }
\end{abstract}



\renewcommand{\theequation}{\thesection.\arabic{equation}}
\renewcommand{\thefigure}{\thesection.\arabic{figure}}

\setcounter{equation}{0}
\section{Introduction and main results}
In this paper, we study the multiplicity and stability of closed characteristics on
any compact convex hypersurface with some kind of symmetries in ${\bf R}^{2n}$.
Let $\Sigma$ be a $C^2$ compact hypersurface in ${\bf R}^{2n}$,
bounding a strictly convex compact set $U$ with non-empty interior, where $n\geq2$.
Denote the set of all such hypersurfaces by $\mathcal{H}(2n)$.
Without loss of generality, we suppose $U$ contains the origin. Let $P$ be a $2n\times2n$ symplectic orthogonal
matrix and $P^k=I_{2n}$, where $k\geq2$.
We denote by $\mathcal {H}_P(2n)$ the set of all P-cyclic symmetric hypersurfaces
in $\mathcal {H}(2n)$, where $\Sigma$ is called {\it P-cyclic symmetric}
if $P\Sigma=\Sigma$, i.e., $x\in \Sigma$ implies $Px\in\Sigma$.
We consider closed characteristics $(\tau,y)$ on $\Sigma$, which are
solutions of the following problem \be
\left\{\matrix{\dot{y}(t)=JN_\Sigma(y(t)),~y(t)\in \Sg,~\forall~t\in
{\bf R}, \cr
               y(\tau)=y(0), ~~~~~~~~~~~~~~~~~~~~~~~~~~~~~~~~~\cr }\right. \lb{1.1}\ee
where $J=\left(
              \begin{array}{cc}
                0 & -I_{n} \\
                I_{n} & 0 \\
              \end{array}
            \right)
$, $I_{n}$ is the identity matrix in ${\bf R}^n$ and
$\mathit{N}_\Sigma(y)$ is the outward normal unit vector  of
$\Sigma$ at $y$ normalized by the condition
$\mathit{N}_\Sigma(y)\cdot y=1$. Here $a\cdot b$ denotes the
standard inner product of $a, b \in {\bf R}^{2n}$. A closed
characteristic $(\tau, y)$ is {\it prime} if $\tau$ is the minimal
period of $y$. Two closed characteristics $(\tau,x)$ and
$(\sigma,y)$ are {\it geometrically distinct}, if $x({\bf
R})\not=y({\bf R})$. We denote by $\mathcal{T}(\Sigma)$ the set of all geometrically distinct closed characteristics
$(\tau,y)$ on $\Sigma$ with $\tau$ being the minimal period of $y$.
A closed characteristic $(\tau, y)$ on $\Sigma \in
\mathcal {H}_P(2n)$ is called {\it P-cyclic symmetric} if $y({\bf
R})=Py({\bf R})$, cf. Proposition 1 of \cite{Zha1}. In this paper, we further assume $ker(P^l-I_{2n})=0$ holds for any $1\leq l< k$.

Let $j: {\bf R}^{2n}\to {\bf R}$ be the gauge function of $\Sigma$,
i.e., $j(\lambda x)=\lambda$ for $x\in\Sigma$ and $\lambda\geq 0$,
then $j\in C^{2}({\bf R}^{2n}\setminus\{0\},{\bf R})\cap C^{1}({\bf
R}^{2n},{\bf R})$ and $\Sigma=j^{-1}(1)$. Fix a constant $\alpha
\in(1, 2)$ and define the Hamiltonian $H_\alpha: {\bf
R}^{2n}\to[0,+\infty)$ by\[ H_\alpha(x) :=j(x)^{\alpha}\] Then
$H_\alpha\in C^{2}({\bf R}^{2n}\setminus\{0\},{\bf R})\cap
C^{0}({\bf R}^{2n},{\bf R})$ is convex and
$\Sigma=H_\alpha^{-1}(1)$. It is well known that the problem
$(1.1)$ is equivalent to the following given energy problem of the
Hamiltonian system
\bea\left\{\begin{array}{ll}
\dot{y}(t)=JH_\alpha^\prime(y(t)), H_\alpha(y(t))=1,~\forall~t\in {\bf R},\\
y(\tau)=y(0).
\end{array}\right. \lb{1.2}\eea
Denote by $\mathcal{T}(\Sigma,\alpha)$ the set of all geometrically distinct solutions
$(\tau,y)$ of the problem $(1.2)$, where $\tau$ is the minimal
period of $y$. Note that elements
in $\mathcal{T}(\Sigma)$ and $\mathcal{T}(\Sigma,\alpha)$ are in one
to one correspondence with each other. Let $(\tau,y)
\in\mathcal{T}(\Sigma,\alpha)$. We call the fundamental solution
$\gamma_y: [0,\tau]\to Sp(2n)$ with $\gamma_y(0)=I_{2n}$ of the
linearized Hamiltonian system\bea
\dot{z}(t)=JH_\alpha^{\prime\prime}(y(t))z(t),~\forall~t \in {\bf
R}.\lb{1.3}\eea the {\it associated symplectic path} of $(\tau,y)$. The
eigenvalue of $\gamma_y(\tau)$ are called  {\it Floquet multipliers}
of $(\tau,y)$. By Proposition 1.6.13 of \cite{Eke1}, the Floquet
multipliers with their multiplicities and Krein type numbers of
$(\tau,y)\in\mathcal{T}(\Sigma,\alpha)$ do not depend on the
particular choice of the Hamiltonian function in $(1.2)$. As in
Chapter 15 of \cite{Lon4}, for any symplectic matrix M, we define the {\it elliptic height} $e(M)$ of
M by the total algebraic multiplicity of all eigenvalues of M on the unit circle
${\bf U}$ in the complex plane ${\bf C}$. And for any $(\tau,y)\in\mathcal{T}(\Sigma,\alpha)$
we define $e(\tau,y)=e(\gamma_y(\tau))$, and call $(\tau,y)$ {\it hyperbolic} if $e(\tau,y)=2$.

There is a long standing conjecture on the number of closed
characteristics on compact convex hypersurfaces in ${\bf
R}^{2n}$(cf.,Page 235 of \cite{Eke1})\bea \,^\#\mathcal{T}(\Sg)\geq n,~~~\forall
~\Sigma\in\mathcal{H}(2n).\eea Since the pioneering works \cite{Rab1}
of P. Rabinowitz and \cite{Wei1} of A. Weinstein in 1978 on the
existence of at least one closed characteristic on every
hypersurface in $\mathcal {H}(2n)$, the existence of multiple closed
characteristics on $\Sigma\in\mathcal{H}(2n)$ has been deeply
studied by many mathematicians. When $n \geq 2$, besides many
results under pinching conditions, in 1987-1988 I. Ekeland-L.
Lassoued, I. Ekeland-H. Hofer, and A. Szulkin (cf. \cite{EkL1},
\cite{EkH1}, \cite{Szu1}) proved \bea \,^\#\mathcal{T}(\Sg)\geq
2,~~~\forall ~\Sigma\in\mathcal{H}(2n).\nn\eea
In \cite{HWZ1} of 1998, H. Hofer-K. Wysocki-E. Zehnder
proved that $\,^{\#}{\T}(\Sg)=2$ or $\infty$ holds for every
$\Sg\in\H(4)$.
In \cite{LoZ1} of 2002,
Y. Long and C. Zhu proved \bea \,^\#\mathcal{T}(\Sg)\geq
[\frac{n}{2}]+1,~~~\forall~ \Sigma\in\mathcal{H}(2n).\nn\eea  In
\cite{WHL1}, W. Wang, X. Hu and Y. Long proved the conjecture (1.4) for $n = 3$.
In \cite{Wan2}, W. Wang proved the conjecture (1.4) for $n = 4$. In
\cite{LLZ1}, C. Liu, Y. Long and C. Zhu proved the conjecture (1.4) when
$\Sg\in\mathcal{H}_P(2n)$ for $P=-I_{2n}$.

In this paper, we prove the conjecture (1.4) for compact convex P-cyclic symmetric hypersurfaces
with $P$ being very general.

{\bf Theorem 1.1.} {\it For every $\Sigma\in\mathcal {H}_P(2n)$, we have
$^\#\mathcal{T}(\Sg)\geq n$.}

Based on the proof of Theorem 1.1, we can further obtain the following results:

{\bf Theorem 1.2.} {\it For every $\Sigma\in\mathcal {H}_P(2n)$ satisfying $^\#\mathcal{T}(\Sg)<+\infty$, there exist at least $2[\frac{n}{2}]$
geometrically distinct non-hyperbolic closed characteristics on $\Sigma$.}

{\bf Theorem 1.3.} {\it For every $\Sigma\in\mathcal {H}_P(2n)$ satisfying $^\#\mathcal{T}(\Sg)=n$ and $k\geq3$,
all of the closed characteristics on $\Sigma$
are P-cyclic symmetric, where $k$ satisfies $P^k=I_{2n}$.}

{\bf Remark 1.4.}
(i) Let $P=-I_{2n}$, our Theorem 1.1 is the same as Theorem 1.1 of \cite{LLZ1}, thus our theorem
extends the main result of \cite{LLZ1} to compact convex hypersurfaces with general symmetries. Our Theorem 1.2
extends and covers Theorem 1.1 of \cite{LLW} which shows Theorem 1.2 holds for $P=-I_{2n}$.
Our Theorem 1.3 is related to the main results of \cite{Wan1} and \cite{LLWZ}, in which the symmetries of closed
characteristics were considered for the special cases $2\leq n\leq 4$.
For more studies about closed characteristics on compact convex P-cyclic symmetric hypersurfaces,
one can also refer to \cite{DoL1, DoL2, LW, Liu1, LWZ, LiZ1, LiZ2, LiZh, Zha1}.

(ii) Our proofs are more complicate than those of \cite{LLZ1} and \cite{LLW}
because we always need to compute Maslov $(P,\omega)$-index for P-cyclic symmetric closed characteristics, which
are not considered elsewhere to study the problems of closed characteristics, thus our proofs are different from those of other papers.

(iii) For the special case $n=2$, one can even obtain the existence of  $\Z_k$-symmetric unknotted periodic orbit on $\Sigma$ which is the binding of an open book decomposition and each page of the open book is a disk-like global surfaces of section when $P=e^{\frac{2\pi}{k}J}$, which has interesting
applications to celestial mechanics, for example, in \cite{Kim} and \cite{Sch}, the H$\acute{e}$non-Heiles Hamiltonian energy level presents $\Z_3$-symmetry and
the Hamiltonian energy level presents $\Z_4$-symmetry in Hill's lunar problem.

This paper is arranged as follows. In Section 2,
we recall briefly the index theory for symplectic paths, especially the Maslov $(P,\omega)$-index theory for symplectic paths and an important
formula of Maslov $(P,\omega)$-index for the associated symplectic paths with convex Hamiltonian systems.
In Section 3, we briefly review a variational structure for closed characteristics on
compact convex hypersurfaces. In Section 4, we use Maslov $(P,\omega)$-index theory to give a key estimation for P-cyclic symmetric closed characteristics in Proposition 4.3 and then
prove our main results.

In this paper, let $\N$, $\Z$, $\Q$, $\R$, and $\C$ denote
the sets of natural integers, integers, rational
numbers, real numbers, and complex numbers respectively.
Denote by $a\cdot b$ and $|a|$ the standard inner product and norm in
$\R^{2n}$. Denote by $\langle\cdot,\cdot\rangle$ and $\|\cdot\|$
the standard $L^2$-inner product and $L^2$-norm. We denote by
$\,^{\#}A$ the number of elements in the set $A$ when it is finite.
We define the functions $[a] = \max{\{k \in {\bf Z} \mid k \leq
a\}}$.

\setcounter{equation}{0}
\section{Index theory for symplectic paths}
In this section, we recall briefly the index theory for symplectic paths
which will be useful in the studies of closed characteristics.

As usual, the symplectic group $\Sp(2n)$ is defined by
$$ \Sp(2n) = \{M\in {\rm GL}(2n,\R)\,|\,M^TJM=J\}, $$
whose topology is induced from that of $\R^{4n^2}$. For $\tau>0$ we are interested
in paths in $\Sp(2n)$:
$$ \P_{\tau}(2n) = \{\ga\in C([0,\tau],\Sp(2n))\,|\,\ga(0)=I_{2n}\}.$$
We consider this path-space equipped
with the $C^0$-topology.
For any $\om\in\U$ the following codimension $1$ hypersurface in $\Sp(2n)$ is
defined in \cite{Lon2}:
$$ \Sp(2n)_{\om}^0 = \{M\in\Sp(2n)\,|\, \det(M-\om I_{2n}))=0\}.  $$
For any $M\in \Sp(2n)_{\om}^0$, we define a co-orientation of $\Sp(2n)_{\om}^0$
at $M$ by the positive direction $\frac{d}{dt}Me^{t J}|_{t=0}$.  Let
\bea
\Sp(2n)_{\om}^{\ast} &=& \Sp(2n)\bs \Sp(2n)_{\om}^0,   \nn\\
\P_{\tau,\om}^{\ast}(2n) &=&
      \{\ga\in\P_{\tau}(2n)\,|\,\ga(\tau)\in\Sp(2n)_{\om}^{\ast}\}, \nn\\
\P_{\tau,\om}^0(2n) &=& \P_{\tau}(2n)\bs  \P_{\tau,\om}^{\ast}(2n).  \nn\eea
For any two continuous arcs $\xi$ and $\eta:[0,\tau]\to\Sp(2n)$ with
$\xi(\tau)=\eta(0)$, their concatenation is defined
as usual by
$$ \eta\ast\xi(t) = \left\{\matrix{
            \xi(2t), & \quad {\rm if}\;0\le t\le \tau/2, \cr
            \eta(2t-\tau), & \quad {\rm if}\; \tau/2\le t\le \tau. \cr}\right. $$
Given any two $2m_k\times 2m_k$ matrices of square block form
$M_k=\left(\matrix{A_k&B_k\cr
                                C_k&D_k\cr}\right)$ with $k=1, 2$,
as in \cite{Lon4}, the $\;\dm$-product of $M_1$ and $M_2$ is defined by
the following $2(m_1+m_2)\times 2(m_1+m_2)$ matrix $M_1\dm M_2$:
$$ M_1\dm M_2=\left(\matrix{A_1&  0&B_1&  0\cr
                               0&A_2&  0&B_2\cr
                             C_1&  0&D_1&  0\cr
                               0&C_2&  0&D_2\cr}\right). \nn$$  
Denote by $M^{\dm k}$ the $k$-fold $\dm$-product $M\dm\cdots\dm M$. Note
that the $\dm$-product of any two symplectic matrices is symplectic. For any two
paths $\ga_j\in\P_{\tau}(2n_j)$ with $j=0$ and $1$, let
$\ga_0\dm\ga_1(t)= \ga_0(t)\dm\ga_1(t)$ for all $t\in [0,\tau]$.

A special path $\xi_n$ is defined by
\bea \xi_n(t) = \left(\matrix{2-\frac{t}{\tau} & 0 \cr
                                             0 &  (2-\frac{t}{\tau})^{-1}\cr}\right)^{\dm n}
        \qquad {\rm for}\;0\le t\le \tau.  \nn\eea

{\bf Definition 2.1.} (cf. \cite{Lon2}, \cite{Lon4}) {\it For any $\om\in\U$ and
$M\in \Sp(2n)$, define
\bea  \nu_{\om}(M)=\dim_{\C}\ker_{\C}(M - \om I_{2n}).  \nn\eea
For any $\tau>0$ and $\ga\in \P_{\tau}(2n)$, define
\bea  \nu_{\om}(\ga)= \nu_{\om}(\ga(\tau)).  \nn\eea

If $\ga\in\P_{\tau,\om}^{\ast}(2n)$, define
\bea i_{\om}(\ga) = [\Sp(2n)_{\om}^0: \ga\ast\xi_n],  \lb{2.1}\eea
where the right hand side of (\ref{2.1}) is the usual homotopy intersection
number, and the orientation of $\ga\ast\xi_n$ is its positive time direction under
homotopy with fixed end points.

If $\ga\in\P_{\tau,\om}^0(2n)$, we let $\mathcal{F}(\ga)$
be the set of all open neighborhoods of $\ga$ in $\P_{\tau}(2n)$, and define
\bea i_{\om}(\ga) = \sup_{U\in\mathcal{F}(\ga)}\inf\{i_{\om}(\beta)\,|\,
                       \beta\in U\cap\P_{\tau,\om}^{\ast}(2n)\}.
               \nn\eea
Then
$$ (i_{\om}(\ga), \nu_{\om}(\ga)) \in \Z\times \{0,1,\ldots,2n\}, $$
is called the index function of $\ga$ at $\om$. }

Note that when $\om=1$, this index theory was introduced by
C. Conley-E. Zehnder in \cite{CoZ1} for the non-degenerate case with $n\ge 2$,
Y. Long-E. Zehnder in \cite{LZe1} for the non-degenerate case with $n=1$,
and Y. Long in \cite{Lon1} and C. Viterbo in \cite{Vit1} independently for
the degenerate case. The case for general $\om\in\U$ was defined by Y. Long
in \cite{Lon2} in order to study the index iteration theory (cf. \cite{Lon4}
for more details and references).

For any symplectic path $\ga\in\P_{\tau}(2n)$ and $m\in\N$,  we
define its $m$-th iteration $\ga^m:[0,m\tau]\to\Sp(2n)$ by
\bea \ga^m(t) = \ga(t-j\tau)\ga(\tau)^j, \qquad
  {\rm for}\quad j\tau\leq t\leq (j+1)\tau,\;j=0,1,\ldots,m-1.
     \nn\eea
We still denote the extended path on $[0,+\infty)$ by $\ga$.

{\bf Definition 2.2.} (cf. \cite{Lon2}, \cite{Lon4}) {\it For any $\ga\in\P_{\tau}(2n)$,
we define
\bea (i(\ga,m), \nu(\ga,m)) = (i_1(\ga^m), \nu_1(\ga^m)), \qquad \forall m\in\N.
   \nn\eea
The mean index $\hat{i}(\ga,m)$ per $m\tau$ for $m\in\N$ is defined by
\bea \hat{i}(\ga,m) = \lim_{k\to +\infty}\frac{i(\ga,mk)}{k}. \lb{2.2}\eea

For any $M\in\Sp(2n)$ and $\om\in\U$, the {\it splitting numbers} $S_M^{\pm}(\om)$
of $M$ at $\om$ are defined by
\bea S_M^{\pm}(\om)
     = \lim_{\ep\to 0^+}i_{\om\exp(\pm\sqrt{-1}\ep)}(\ga) - i_{\om}(\ga),
   \lb{2.3}\eea
for any path $\ga\in\P_{\tau}(2n)$ satisfying $\ga(\tau)=M$.}

Based on the index theory above, the Maslov $(P,\omega)$-index was defined as follows:

{\bf Definition 2.3.} (cf. \cite{LT1}) {\it For any $P\in Sp(2n)$, $\omega\in \U$ and $\gamma \in \mathcal {P}_{\tau} (2n)$, the Maslov
$(P, \omega)$-index is defined by}\bea i_\omega^P(\gamma)=i_\omega(P^{-1}\gamma*\xi)-i_\omega(\xi),\nn\eea
where $\xi\in \mathcal {P}_{\tau} (2n)$ such that $\xi(\tau)=P^{-1}\gamma(0)=P^{-1}$, and $(P, \omega)$-nullity $\nu^P_
{\omega}(\gamma)$ is defined by \bea \nu_\omega^P(\gamma)=dim_{{\bf C}} ker_{{\bf
C}}(\gamma(\tau) - \omega P).\nn\eea

For any $M \in Sp(2n)$ and $\omega \in {\bf U}$, {\it the splitting
numbers} $_PS^{\pm}_M (\omega)$ of $M$ at $(P,\omega)$ are defined in Definition 2.4 of \cite{LT2} as follows
\begin{eqnarray}
_PS^{\pm}_M (\omega)=\lim_{\epsilon\rightarrow 0^+}{i^P_{\omega
\exp{(\pm\sqrt{-1}\epsilon)}}(\gamma)-i^P_{\omega}(\gamma)},
\lb{2.4}\end{eqnarray}
for any path $\gamma \in \mathcal {P}_{\tau} (2n)$ satisfying
$\gamma(\tau) = M$.

Note that the Maslov P-index theory
for a symplectic path was first studied by Y. Dong and C. Liu in \cite{Dong, LiuC} independently for any symplectic matrix P with different
treatment. The Maslov P-index theory was generalized in \cite{LT1} to the Maslov
$(P, \omega)$-index theory for any $P\in Sp(2n)$ and all $\omega\in \U$. The iteration theory of
$(P, \omega)$-index theory was studied in \cite{LT2}. When $\omega= 1$, the Maslov $(P, \omega)$-index theory
coincides with the Maslov P-index theory. In order to study the properties of P-cyclic symmetric closed characteristics, we will use Maslov $(P,\omega)$-index theory for symplectic paths in Section 4.

Let $\Omega^0(M)$ be the path connected component containing $M =
\gamma(\tau)$ of the set\begin{eqnarray} \Omega(M) = \{N \in Sp(2n)
\mid  \sigma(N)
\cap {\bf U} = \sigma(M) \cap {\bf U}~and~~~~\nn\\
\nu_{\lambda}(N)=\nu_{\lambda}(M),\forall \lambda \in \sigma(M) \cap
{\bf U}\}\nn\end{eqnarray} Here $\Omega^0(M)$ is called the {\it
homotopy component} of $M$ in $Sp(2n)$.

In \cite{Lon2}-\cite{Lon4}, the following symplectic matrices were introduced as
basic normal forms: \begin{eqnarray}D(\lambda)=\left(
                                                    \begin{array}{cc}
                                                      \lambda & 0 \\
                                                      0 & \lambda^{-1} \\
                                                    \end{array}
                                                  \right),
                                                  \lambda=\pm 2,\lb{2.5}\\
N_1(\lambda,b)=\left(
 \begin{array}{cc}
  \lambda & b \\
   0& \lambda \\
   \end{array}
   \right),\lambda=\pm1,b=\pm1,0,\lb{2.6}\\
R(\theta)=\left(
            \begin{array}{cc}
              \cos{\theta} & -\sin{\theta} \\
              \sin{\theta} & \cos{\theta} \\
            \end{array}
          \right),\theta\in (0,\pi)\cup(\pi,2\pi),\lb{2.7}\\
N_2(\omega,B)=\left(
                \begin{array}{cc}
                  R(\theta) & B \\
                  0 &  R(\theta) \\
                \end{array}
              \right),\theta\in (0,\pi)\cup(\pi,2\pi),
\lb{2.8}\end{eqnarray}
where $B=\left(
           \begin{array}{cc}
             b_1 & b_2 \\
             b_3 & b_4 \\
           \end{array}
         \right)
$ with $b_i\in {\bf R}$ and $b_2\neq b_3$, $\omega=e^{\theta\sqrt{-1}}$.

{\bf Lemma 2.4.} (cf. Theorem 1.8.10 of \cite{Lon4}) {\it For any $M
\in Sp(2n)$, there is a path $f : [0, 1] \rightarrow \Omega^0(M)$
such that $f(0) = M$ and\begin{eqnarray} f(1)=M_1\diamond \cdots
\diamond M_l,\nn
\end{eqnarray}
where each $M_i$ is a basic normal form listed in (2.5)-(2.8) for
$1 \leq i \leq l$.}

Splitting numbers possess the following properties:

{\bf Lemma 2.5.} (cf. Lemma 9.1.5, Lemma 9.1.6, List 9.1.12 and Corollary 9.2.4 of \cite{Lon4})
{\it Splitting numbers $S_M^{\pm}(\om)$ are well defined, i.e., they are independent of the choice
of the path $\ga\in\P_\tau(2n)$ satisfying $\ga(\tau)=M$ appeared in (\ref{2.3}). For $M\in\Sp(2n)$, splitting numbers
$S^{\pm}_N(\omega)$ are constant for all $N \in \Omega^0(M)$. Moreover,
there hold
\begin{eqnarray}
S_M^{+}(\om) &=& S_M^{-}(\bar{\om}),~\forall~\omega
\in {\bf U}.\nn\\
S_M^{\pm}(\om) &=& 0, \qquad {\it if}\;\;\om\not\in\sg(M).  \nn\\
S_{M^m}^{\pm}(z) &=& \sum_{\omega^m=z}S_M^{\pm}(\omega),~\forall~z
\in {\bf U}, m\in\N.\nn\\
S_{N_1(1,a)}^+(1) &=& \left\{\matrix{1, &\quad {\rm if}\;\; a\ge 0, \cr
0, &\quad {\rm if}\;\; a< 0. \cr}\right. \nn\\S_{N_1(-1,a)}^+(-1) &=& \left\{\matrix{1, &\quad {\rm if}\;\; a\leq 0, \cr
0, &\quad {\rm if}\;\; a> 0. \cr}\right. \nn\\
(S_{R(\theta)}^+(e^{\sqrt{-1}\theta}), S_{R(\theta)}^-(e^{\sqrt{-1}\theta}))&=&(0,1), \theta \in (0, \pi) \cup (\pi, 2\pi) \nn
\\
S_{N_2(e^{\sqrt{-1}\theta},B)}^{\pm}(e^{\sqrt{-1}\theta}) &=& \left\{\matrix{1, &\quad {\rm if}\;\; (b_2-b_3)\sin{\theta}>0, \cr
0, &\quad {\rm if}\;\; (b_2-b_3)\sin{\theta}<0, \cr}\right. \nn\eea
where $B=\left(
           \begin{array}{cc}
             b_1 & b_2 \\
             b_3 & b_4 \\
           \end{array}
         \right)
$ with $b_i\in {\bf R}$ and $b_2\neq b_3$.
For any $M_i\in\Sp(2n_i)$ with $i=0$ and $1$, there holds }
\bea S^{\pm}_{M_0\dm M_1}(\om) = S^{\pm}_{M_0}(\om) + S^{\pm}_{M_1}(\om),
    \qquad \forall\;\om\in\U. \nn\eea

{\bf Lemma 2.6.}(cf. Lemma 2.5 and Lemma 2.6 of \cite{LT2}) {\it For any $M\in Sp(2n)$ and $\omega\in {\bf U}$, the
splitting numbers $_PS^{\pm}_M (\omega)$ are well defined, i.e., they are independent of the choice
of the path $\ga\in\P_\tau(2n)$ satisfying $\ga(\tau)=M$ appeared in (\ref{2.4}). And
the following properties hold:

(i) $_PS^{\pm}_M (\omega)=S^{\pm}_{P^{-1}M}(\omega)-S^{\pm}_{P^{-1}}(\omega)$.

(ii) $_PS^{+}_M (\omega)= {}_PS^{-}_M (\bar{\omega})$.

(iii) $_PS^{\pm}_M (\omega)= {}_PS^{\pm}_N (\omega)$ if $P^{-1}N\in
\Omega^0(P^{-1}M)$.

(iv) $_PS^{\pm}_{M_1\diamond
M_2}(\omega)= {}_{P_1}S^{\pm}_{M_1}(\omega)+ {}_{P_2}S^{\pm}_{M_2}(\omega)$
for $M_j, P_j\in Sp(2n_j)$ with $n_j\in\{1,\cdots,n\}$ satisfying
$P=P_1\diamond P_2$ and $n_1+n_2=n$.

(v) $_PS^{\pm}_M (\omega)=0$ if $\omega\notin \sigma(P^{-1}M)\cup\sigma(P^{-1})$.}

For any symplectic path $\ga\in\P_{\tau}(2n)$ and $m\in\N$,  we define the $m$-times iteration path $\ga_P^m:[0,m\tau]\to\Sp(2n)$
of $\gamma$ by
\bea \ga_P^m(t) =\left\{\matrix{
            \gamma(t), & \quad t\in[0,\tau], \cr
            P\gamma(t-\tau)P^{-1}\gamma(\tau), &  \quad t\in[\tau, 2\tau], \cr
            P^2\gamma(t-2\tau)(P^{-1}\gamma(\tau))^2, &  \quad t\in[2\tau, 3\tau], \cr
            \cdots\cr
            P^{m-1}\gamma(t-(m-1)\tau)(P^{-1}\gamma(\tau))^{m-1}, &  \quad t\in[(m-1)\tau, m\tau].\cr}\right. \eea

By Lemma 2.7 of \cite{LT2}, we have the Bott-type formula for Maslov $(P,\omega)$-index.:

{\bf Lemma 2.7.} {\it For any $\gamma \in \mathcal {P}_{\tau} (2n)$, $z \in {\bf U}$
and $m\in\N$, we have\bea i_z^{P^m}(\gamma_P^m)&=&\sum_{\omega^{m}=z}
i^P_{\omega}(\gamma),\nn\\ \nu_z^{P^m}(\gamma_P^m)&=&\sum_{\omega^{m}=z}
\nu^P_{\omega}(\gamma).\nn\eea}
The following formula of Maslov $(P,\omega)$-index for the associated symplectic paths with convex Hamiltonian systems will play an important role
in the proof of Proposition 4.3 below.

{\bf Lemma 2.8.}(cf. Lemma 3.5 of \cite{LWZ}) {\it Assume $A(t)\in GL(\R^{2n})$ is positive definite for $t\in[0,\tau]$, let $\gamma\equiv\gamma_A\in \mathcal {P}_{\tau} (2n)$
be the fundamental solution of the linearized Hamiltonian system $\dot{y}(t)=JA(t)y(t)$. Then we have
\bea i_\omega^P(\gamma)=\nu_\omega(P^{-1})+\sum_{0<s<\tau}\nu_\omega^P(\gamma(s)), \forall \omega\in\U,\nn\eea
where $P$ satisfies $ker(P-I_{2n})=0$ as assumed in \cite{LWZ}.}

\setcounter{equation}{0}
\section{Variational properties for closed characteristics }

In this section, we describe the variational properties for closed
characteristics.

To solve the given energy problem (1.2), we
consider the fixed period problem
\be \left\{\matrix{\dot{x}(t)=JH_\alpha^\prime(x(t)), \cr
     x(1)=x(0).         \cr }\right. \lb{3.1}\ee

Define
\bea L_0^{\frac{\alpha}{\alpha-1}}(S^1,\R^{2n})
  =\{u\in L^{\frac{\alpha}{\alpha-1}}(S^1,\R^{2n})\,|\,\int_0^1udt=0\}. \nn\eea
The corresponding Clarke-Ekeland dual action functional is defined by
\be \Phi(u)=\int_0^1\left(\frac{1}{2}Ju\cdot Mu+H_\alpha^{\ast}(-Ju)\right)dt,
    \qquad \forall\;u\in L_0^{\frac{\alpha}{\alpha-1}}(S^1,\R^{2n}), \lb{3.2}\ee
where $Mu$ is defined by $\frac{d}{dt}Mu(t)=u(t)$ and $\int_0^1Mu(t)dt=0$,
$H_\alpha^\ast$ is the Fenchel transform of $H_\alpha$ defined by
$H_\alpha^\ast(y)=\sup\{x\cdot y-H_\alpha(x)\;|\; x\in \R^{2n}\}$.
By Theorem 5.2.8 of \cite{Eke1}, $\Phi$ is $C^1$ on
$L_0^{\frac{\alpha}{\alpha-1}}(S^1,\,\R^{2n})$ and satisfies
the Palais-Smale condition.  Suppose
$x$ is a solution of (\ref{3.1}). Then $u=\dot{x}$ is a critical point
of $\Phi$. Conversely, suppose $u$ is a critical point of $\Phi$.
Then there exists a unique $\xi\in\R^{2n}$ such that $x_{u}=Mu-\xi$ is a
solution of (\ref{3.1}). In particular, solutions of (\ref{3.1}) are in
one to one correspondence with critical points of $\Phi$. Moreover,
$\Phi(u)<0$ for every critical point $u\not= 0$ of $\Phi$. In addition, we have a natural $S^1$-action on $L_0^{\frac{\alpha}{\alpha-1}}(S^1,\; \R^{2n})$ defined by
$\th\cdot u(t)=u(\th+t)$ for all $\th\in S^1$ and $t\in\R$. Clearly
$\Phi$ is $S^1$-invariant.

Let $h=H_{\alpha}(x_{u}(t))$ and $1/m$ be
the minimal period of $x_{u}$ for some $m\in {\bf N}$. Define
\begin{eqnarray}
y_u(t)=h^{\frac{-1}{\alpha}}x_u(h^{\frac{2-\alpha}{\alpha}}t) ~~~
and ~~~\tau=\frac{1}{m}h^{\frac{\alpha-2}{\alpha}}.\nn
\end{eqnarray}
Then there hold $y_{u}(t)\in \Sigma$ for all $t\in {\bf R}$ and
$(\tau,y_{u})\in \mathcal{T}(\Sigma,\alpha)$. Note that the period 1
of $x_{u}$ corresponds to the period $m\tau$ of the solution
$(m\tau,y_{u}^{m})$ of (1.2) with minimal period $\tau$.
On the other hand, every solution $(\tau,y)\in
\mathcal{T}(\Sigma,\alpha)$ gives rise to a sequence
$\{x_{y}^{m}\}_{m\in {\bf N}}$ of solutions of the problem $(3.1)$,
and a sequence $\{u_{y}^{m}\}_{m\in {\bf N}}$ of critical points of
$\Phi$ defined by
\begin{eqnarray}
x_y^ m(t)=(m\tau)^{\frac{-1}{2-\alpha}}y(m\tau t)\nn\\
u_y^ m(t)=(m\tau)^{\frac{1-\alpha}{2-\alpha}}\dot {y}(m\tau t)
\end{eqnarray}

Suppose $u$ is a nonzero critical point of $\Phi$. Then
the formal Hessian of $\Phi$ at $u$ on $L_0^2(S^1,\R^{2n})$ is defined by
\bea Q(v,\; v)=\int_0^1 (Jv\cdot Mv+(H_\alpha^\ast)^{\prime\prime}(-Ju)Jv\cdot Jv)dt,\eea
which defines an orthogonal splitting $L_0^2(S^1,\R^{2n})=E_-\oplus E_0\oplus E_+$
 into negative, zero and positive subspaces. The Ekeland
index of $u$ is defined by $i(u)=\dim E_-$ and the nullity of $u$ is
defined by $\nu(u)=\dim E_0$. The inequality $1\le \nu(u)\le 2n$
always holds, cf. P.219 of  \cite{Eke1}.

For a closed characteristic $(\tau,y)\in \mathcal{T}(\Sigma,\alpha)$, we denote by
$y^m\equiv (m\tau, y)$ the $m$-th iteration of $y$ for $m\in\N$.
Then we define the index $i(y^m)$ and nullity $\nu(y^m)$
of $(m\tau,y)$ for $m\in\N$ by
\bea i(y^m)=i(u_y^m), \qquad \nu(y^m)=\nu(u_y^m). \eea
The mean index of $(\tau,y)$ is defined by
\be \hat{i}(y)=\lim_{m\rightarrow\infty}\frac{i(y^m)}{m}. \ee
Note that by Corollary 8.3.2 and Lemma 15.3.2
of \cite{Lon4}, there always holds
\bea \hat{i}(y)>2.\eea
We define via Definition 2.2 the following
\bea  S^+(y) &=& S_{\ga_y(\tau)}^+(1),  \\
  (i(y,m), \nu(y,m)) &=& (i(\ga_y,m), \nu(\ga_y,m)), \\
   \hat{i}(y,m) &=& \hat{i}(\ga_y,m), \eea
for all $m\in\N$, where $\ga_y$ is the associated symplectic path of $(\tau,y)$, i.e., the
fundamental solution of the linearized Hamiltonian system of (\ref{1.3}) at $(\tau,y)$.

{\bf Lemma 3.1.} (cf. Lemma 1.1 of \cite{LoZ1}, Theorem 15.1.1 of \cite{Lon4}) {\it Suppose
$(\tau,y)\in \T(\Sigma, \alpha)$. Then we have
\bea i(y^m)\equiv i(m\tau ,y)=i(y, m)-n,\quad \nu(y^m)\equiv\nu(m\tau, y)=\nu(y, m),
       \qquad \forall m\in\N.\nn\eea
In particular, (\ref{2.2}), (3.6) and (3.10)
coincide, thus we simply denote them by $\hat i(y)$.}

Since the Ekeland index is a Morse-type index which is non-negative, then by Lemma 3.1 we have

{\bf Corollary 3.2.}(cf. Corollary 15.1.4 of \cite{Lon4}) It holds
that
\begin{eqnarray}
i(y, 1)\geq n,~\forall~(\tau,y)\in \mathcal{T}(\Sigma,\alpha).\nn
\end{eqnarray}

By Corollary 3.1 of \cite{LoZ1}, we have the monotonicity for the index iterations of closed characteristics:

{\bf Lemma 3.3.}  {\it For any
$(\tau,y)\in \T(\Sigma, \alpha)$, $m\in {\bf N}$, we have}
\bea i(y, m+1)-i(y, m)&\geq& 2,\\ i(y, m+1)+\nu(y, m+1)-1&\geq& i(y, m+1)> i(y, m)+\nu(y, m)-1.
       \eea

Following Section V.3 of \cite{Eke1}, denote by ``ind'' the
$S^1$-action cohomology index theory for $S^1$-invariant subset of
$L_0^2(S^1,\R^{2n})$ defined in \cite{Eke1}. For $[\Phi]_c\equiv\{u\in
L_0^2(S^1,\R^{2n})\mid \Phi(u)\leq c\}$ define
\begin{eqnarray}
c_k=\inf{\{c<0\mid ind([\Phi]_c)\geq k\}}.\nn
\end{eqnarray}
Then by Proposition 3
in P.218 of \cite{Eke1}, we have

{\bf Lemma 3.4.} {\it Every $c_i$ is a critical value of $\Phi$. If
$c_i=c_j$ for some $i<j$, then there are infinitely many geometrically
distinct closed characteristics on $\Sg$.}

By Theorem 4 in P.219 of \cite{Eke1}, we have the following

{\bf Lemma 3.5.} {\it For any
given $k\in {\bf N}$, there exists $(\tau,y)\in
\mathcal{T}(\Sigma,\alpha)$ and $m\in {\bf N}$ such that for $u_y^m$
defined by $(3.3)$}
\begin{eqnarray}
\Phi^\prime( u_y^m)=0~~~~~~and~~~~~~\Phi(
u_y^m)=c_k,\nn\\i(u_y^m)\leq 2k-2\leq i(u_y^m)+\nu(u_y^m)-1.\nn
\end{eqnarray}

Combining Lemma 3.1, Lemma 3.4 and Lemma 3.5, we obtain

{\bf Lemma 3.6.}(cf. Lemma 3.1 of  \cite{LoZ1}) {\it Suppose  $^\#\mathcal{T}(\Sg)<+\infty$, there exists an
injection map $p = p(\Sigma,\alpha) :\N\rightarrow \mathcal{T}(\Sg, \alpha)\times\N$ such that for  any
$k\in {\bf N}$, $(\tau,y)\in
\mathcal{T}(\Sigma,\alpha)$ and $m\in\N$ satisfying $p(k) = ((\tau, y),m)$, there hold
\begin{eqnarray}
\Phi^\prime( u_y^m)=0~~~~~~and~~~~~~\Phi(
u_y^m)=c_k,\nn\\i(y,m)\leq 2k-2+n\leq i(y,m)+\nu(y,m)-1,
\end{eqnarray}
where $u_y^m$ is defined by (3.3).}

\setcounter{equation}{0}
\section{Proof of the main results}
Before we give the proofs of our main results, we first consider the normal forms of the cyclic symplectic orthogonal matrix $P$,
which is crucial for our proof of Proposition 4.3 below.

{\bf Proposition 4.1.} {\it Let $P$ be a symplectic orthogonal matrix satisfying $P^k=I_{2n}$ for some integer $k\geq 2$, then
there exists some $Q\in Sp(2n)$ such that
\bea QPQ^{-1}=R(\theta_1)\diamond R(\theta_2)\diamond\cdots \diamond R(\theta_n),\nn\eea
where $\theta_i\in[0,2\pi)$.}

{\bf Proof.} By Lemma 2.4, there exists some $Q\in Sp(2n)$ such that \bea QPQ^{-1}
&=& N_1(1,1)^{\dm p_-}\,\dm\,N_1(1,-1)^{\dm p_+}
  \dm\,N_1(-1,1)^{\dm q_-}\,\dm\,N_1(-1,-1)^{\dm q_+} \nn\\
&&\dm\,N_2(e^{\aa_{1}\sqrt{-1}},A_{1})\,\dm\,\cdots\,\dm\,N_2(e^{\aa_{r_{\ast}}\sqrt{-1}},A_{r_{\ast}})\nn\\
&&\dm\,R(\th_1)\,\dm\,\cdots\,\dm\,R(\th_r)\dm\,D(\pm 2)^{\dm h},\nn\eea
where $\aa_{j}\in(0,\pi)\cup(\pi,2\pi)$ for $1\le j\le r_{\ast}$, $\th_{j^\prime}\in[0,2\pi)$ for $1\le j^\prime\le r$ and non-negative integers
$p_-, p_+,q_-, q_+,r,r_\ast,h$ satisfy
\bea p_- + p_+ + q_- + q_+ + r + 2r_{\ast}+ h = n. \nn\eea
Note that $I_{2n}=QP^kQ^{-1}=(QPQ^{-1})^k$, then we must have
\bea QPQ^{-1}
&=& N_2(e^{\aa_{1}\sqrt{-1}},A_{1})\,\dm\,\cdots\,\dm\,N_2(e^{\aa_{r_{\ast}}\sqrt{-1}},A_{r_{\ast}})\nn\\
&&\dm\,R(\th_1)\,\dm\,\cdots\,\dm\,R(\th_r),\\
\left(N_2(e^{\aa_{j}\sqrt{-1}},A_{j})\right)^k&=&I_4,\quad 1\leq j\leq r_{\ast}.\eea
On the other hand,
by direct computation we have \bea \left(N_2(e^{\aa_{j}\sqrt{-1}},A_{j})\right)^k&=&\left(
                \begin{array}{cc}
                  R(k\aa_{j}) & B \\
                  0 &  R(k\aa_{j}) \\
                \end{array}
              \right),\\B&=&\sum_{i=0}^k  R(i\aa_{j})A_{j}R((k-i)\aa_{j}).\eea
Comparing (4.2) with (4.3), we have $R(k\aa_{j})=I_2$ and $B=0$ which together with (4.4) implies
\bea 0=R(\aa_{j})BR(-\aa_{j})=\sum_{i=1}^{k+1}  R(i\aa_{j})A_{j}R((k-i)\aa_{j}).\eea
Comparing it with (4.4), we get $A_j=R(\aa_{j})A_{j}R(-\aa_{j})$. Combining it
with (4.5), we obtain $A_j=0$. Thus $r_{\ast}=0$ and by (4.1) we complete the proof.\hfill\hb

In the following, we fix a $\Sg\in\mathcal{H}_P(2n)$ satisfying $P^k=I_{2n}$ and $ker(P^l-I_{2n})=0$ holds for any $1\leq l< k$,
and assume that there exist only finitely many geometrically distinct closed characteristics on $\Sigma$,
i.e., $\mathcal{T}(\Sg) = \{(\tau_j, y_j)\}_{1\le j\le S}$.
We denote by $\ga_j\equiv \gamma_{y_j}$ the associated symplectic path of $(\tau_j,\,y_j)$ on $\Sg$ for
$1\le j\le S$. Then by Lemma 1.3 of \cite{LoZ1} or Lemma 15.2.4 of \cite{Lon4}, there exist $P_j\in \Sp(2n)$
and $M_j\in \Sp(2n-2)$ such that
\bea \ga_j(\tau_j)=P_j^{-1}(N_1(1,\,1)\dm M_j)P_j, \quad\forall\; 1\le j\le S.
   \eea
{\bf Proposition 4.2.}  {\it Suppose $(\tau,\, y)\in\mathcal{T}(\Sigma, \,\alpha)$,
then $(\tau,\, Py)\in\mathcal{T}(\Sigma, \,\alpha)$
and either $\mathcal{O}(y)=\mathcal{O}(Py)$ or
$\mathcal{O}(y)\cap\mathcal{O}(Py)=\emptyset$,
where $\mathcal{O}(P y)=\{P y(t)|\, t\in\R\}$.
Moreover, if $\mathcal{O}(y)\cap\mathcal{O}(Py)\neq\emptyset$,
then we have
\bea  y\left(t+\frac{\tau}{k}\right)=P^ly(t),\qquad\forall t\in\R.\eea
for some $1\leq l<k$.}

{\bf Proof.} Since $\Sigma=P\Sigma$, then
\bea H_\alpha(Py)&=&H_\alpha(y),\\
H_\alpha^\prime(Py)&=&PH_\alpha^\prime(y),\\
H_\alpha^{\prime\prime}(y)&=&P^{T}H_\alpha^{\prime\prime}(Py)P.\eea
Note that $PJ=JP$ since $P$ is a symplectic orthogonal matrix, which together with (1.2) and (4.9) yields
$(\tau,\, Py)\in\mathcal{T}(\Sigma, \,\alpha)$ for any $(\tau,\, y)\in\mathcal{T}(\Sigma, \,\alpha)$.

If $\mathcal{O}(y)\cap\mathcal{O}(Py)\neq\emptyset$, there exists an $s \in [0, \tau)$ such that $y(s) = Py(0)$.
Because $y(s+ t)$ and $Py(t)$ satisfy the same Hamiltonian system
$$ \dot{x}=JH_\alpha^\prime(x),$$
and they have the same initial value for $t=0$,
we have \bea y(t+s) = Py(t),\quad \forall t\in\R, \eea
which implies $\mathcal{O}(y)=\mathcal{O}(Py)$. Combining (4.11) with
the assumption that $P^k=I_{2n}$, we get $y(t+ks) = y(t)$,
then there holds
\bea s=\frac{j}{k}\tau,\eea
for some $0\leq j<k$.

{\bf Claim.} $j$ and $k$ are co-prime.

If otherwise, there exists some  $1\leq l< k$
such that $jl\equiv 0(\mod k)$, then $y(t)=y(t+ls) = P^ly(t)$ by (4.11)-(4.12).
Note that $y(t)\neq0$ for any $t\in\R$ and $ker(P^l-I_{2n})=0$ for $1\leq l< k$, we get a contradiction
and the above claim holds.

By this claim, there exists some  $1\leq l< k$ such that $jl\equiv 1(\mod k)$, hence we obtain $y(t+\frac{\tau}{k})=y(t+ls) = P^ly(t)$ by (4.11)-(4.12).
The proof is complete.\hfill\hb

As in \cite{Zha1}, we call a closed characteristic $(\tau,\, y)$
on $\Sg\in\mathcal{H}_P(2n)$ {\it P-cyclic symmetric}
if $\mathcal{O}(y)=\mathcal{O}(Py)$.
Thus by Proposition 4.2, we obtain that if $(\tau,\, y)$ is not P-cyclic symmetric, then
$(\tau,\, Py)$ and $(\tau,\, y)$ are geometrically distinct.

For $\Sg\in\mathcal{H}_P(2n)$, we have $\Phi(Pu)=\Phi(u)$. In fact, by (4.8), we have $H_\alpha^\ast(Py)=H_\alpha^\ast(y)$
which together with (3.2) yields
\bea \Phi(Pu)&=&\int_0^1\left(\frac{1}{2}JPu\cdot MPu+H_\alpha^{\ast}(-JPu)\right)dt\nn\\
&=&\int_0^1\left(\frac{1}{2}PJu\cdot PMu+H_\alpha^{\ast}(-PJu)\right)dt\nn\\
&=&\int_0^1\left(\frac{1}{2}Ju\cdot Mu+H_\alpha^{\ast}(-Ju)\right)dt\nn\\&=&\Phi(u),
    \qquad \forall\;u\in L_0^{\frac{\alpha}{\alpha-1}}(S^1,\R^{2n}), \eea
where we used the fact that $PJ=JP$.
Note that if $u\in L_0^{\frac{\alpha}{\alpha-1}}(S^1,\R^{2n})$ is the critical point of
$\Phi$ corresponding to a closed characteristic $(\sigma, z)$, then
$Pu\in L_0^{\frac{\alpha}{\alpha-1}}(S^1,\R^{2n})$ is the critical point of
$\Phi$ corresponding to the closed characteristic $(\sigma, Pz)$, thus $u$ and $Pu$ have the same critical values.
Moreover, by (3.4)-(3.5), Lemma 3.1 and (4.10), we have
\bea i(z, 1)=i(Pz,1),\quad \nu(z,1)=\nu(Pz,1).\eea

Now we give a key estimation for P-cyclic symmetric closed characteristics, which is a crucial step for proving our main results.
As far as we know, it is not considered by other papers.

{\bf Proposition 4.3.}  {\it For any P-cyclic symmetric closed characteristic $(\tau,\, y)\in\mathcal{T}(\Sigma, \,\alpha)$
on $\Sg\in\mathcal{H}_P(2n)$, we have
\bea    i(y,\, 1) + 2S^+(y) - \nu(y,\,1) \ge n,\eea
where we use the notations in (3.8)-(3.9).}

{\bf Proof.} By (4.7), we have \bea  y\left(t+\frac{\tau}{k}\right)=P^ly(t),\qquad\forall t\in\R,\eea
for some $1\leq l<k$. Let $\bar{P}=P^l$, then $\Sigma$ is still $\bar{P}$-cyclic symmetric, and (4.16) becomes \bea  y\left(t+\frac{\tau}{k}\right)=\bar{P}y(t),\qquad\forall t\in\R.\eea
In the following, we use $\bar{P}$ instead of $P$ considered in the $(P, \omega)$ index theory of Section 2.
Let $\gamma_y : [0,\,\tau]\rightarrow \Sp(2n)$ with $\gamma_y(0)=I_{2n}$ be the associated symplectic path of $(\tau,y)$ and $\gamma=\gamma_y|_{[0,\,\tau/k]}$.
From (4.10) and (4.17), we obtain $H_\alpha^{\prime\prime}(y(t))=\bar{P}^{T}H_\alpha^{\prime\prime}(y(t+\frac{\tau}{k}))\bar{P}$, then by (1.5)-(1.6) of \cite{LT1} we have
\bea \gamma_y=\gamma_{\bar{P}}^k,\eea
which is the $k$-times iteration path of $\gamma$ defined as in (2.9).

By Proposition 4.1 and the assumption that $ker(P^l-I_{2n})=0$ for $1\leq l< k$, there exists some $Q\in Sp(2n)$ such that
\bea  Q\bar{P}Q^{-1}=R(\theta_1)\diamond R(\theta_2)\diamond\cdots \diamond R(\theta_n),\quad 0<\frac{\theta_i}{\pi}< 2.\eea Noticing that $\bar{P}^k=I_{2n}$,
then by (4.19) we have \bea \sum_{\omega^{k}=1, \omega\neq 1} \nu_\omega({\bar{P}}^{-1})=2n.\eea
Let $\omega_1=e^{\frac{2\pi}{k}\sqrt{-1}}$ and denote the eigenvalues of $\bar{P}^{-1}M$ between 1 and $\omega_1$ lying on the upper
semi-circle in $\U$ by $\alpha_1, \alpha_2,\cdots, \alpha_r$ anticlockwise, where $M=\gamma(\tau/k)$.
Let $\alpha_0=1$ and $\alpha_{r+1}=\omega_1$.
By the definitions of splitting numbers, we have
\bea i^{\bar{P}}_{\alpha_i}(\gamma)+{}_{\bar{P}}S^{+}_M (\alpha_i)=i^{\bar{P}}_{\alpha_{i+1}}(\gamma)+{}_{\bar{P}}S^{-}_M (\alpha_{i+1}),i=0,1,...,r.\eea
Note that by Lemma 2.6(i), it follows that\bea {}_{\bar{P}}S^{\pm}_M (\omega)=S^{\pm}_{{\bar{P}}^{-1}M}(\omega)-S^{\pm}_{{\bar{P}}^{-1}}(\omega).\eea
Then by (4.21)-(4.22) we have
\bea i^{\bar{P}}_{\alpha_{i+1}}(\gamma)-i^{\bar{P}}_{\alpha_i}(\gamma)&=&S^{+}_{{\bar{P}}^{-1}M}(\alpha_i)-S^{-}_{{\bar{P}}^{-1}M}(\alpha_{i+1})\nn\\
&&+S^{-}_{{\bar{P}}^{-1}}(\alpha_{i+1})-S^{+}_{{\bar{P}}^{-1}}(\alpha_i),i=0,1,...,r,\nn\eea
which implies
\bea i^{\bar{P}}_{\omega_1}(\gamma)-i^{\bar{P}}_{1}(\gamma)&=& i^{\bar{P}}_{\alpha_{r+1}}(\gamma)-i^{\bar{P}}_{\alpha_0}(\gamma)\nn\\
&=&\sum_{i=0}^r (i^{\bar{P}}_{\alpha_{i+1}}(\gamma)-i^{\bar{P}}_{\alpha_i}(\gamma))\nn\\
&=&\sum_{i=0}^rS^{+}_{{\bar{P}}^{-1}M}(\alpha_i)-\sum_{i=0}^rS^{-}_{{\bar{P}}^{-1}M}(\alpha_{i+1})\nn\\
&&+\sum_{i=0}^rS^{-}_{{\bar{P}}^{-1}}(\alpha_{i+1})-\sum_{i=0}^rS^{+}_{{\bar{P}}^{-1}}(\alpha_i).\eea
Note that $S^{\pm}_{{\bar{P}}^{-1}}(\alpha_i)=0$ for $0\leq i\leq r$ since $\alpha_i$ is not an eigenvalue of ${\bar{P}}^{-1}$,
then (4.23) becomes\bea i^{\bar{P}}_{1}(\gamma)=i^{\bar{P}}_{\omega_1}(\gamma)+\sum_{i=0}^rS^{-}_{{\bar{P}}^{-1}M}(\alpha_{i+1})-
\sum_{i=0}^rS^{+}_{{\bar{P}}^{-1}M}(\alpha_i)-S^{-}_{{\bar{P}}^{-1}}(\omega_1).\eea
On the other hand, since $ker(\bar{P}-I_{2n})=0$, then it follows from Lemma 2.8 that \bea i_\omega^{\bar{P}}(\gamma)\geq \nu_\omega({\bar{P}}^{-1}).\eea
Combining (4.24) with (4.25), we have
\bea i^{\bar{P}}_{1}(\gamma)&\geq& \sum_{i=0}^rS^{-}_{{\bar{P}}^{-1}M}(\alpha_{i+1})-
\sum_{i=0}^rS^{+}_{{\bar{P}}^{-1}M}(\alpha_i)\nn\\&=&\sum_{i=1}^rS^{-}_{{\bar{P}}^{-1}M}(\alpha_{i})-
\sum_{i=1}^rS^{+}_{{\bar{P}}^{-1}M}(\alpha_i)+S^{-}_{{\bar{P}}^{-1}M}(\omega_1)-S^{+}_{{\bar{P}}^{-1}M}(1),\eea
where we used the fact that $\nu_{\omega_1}({\bar{P}}^{-1})\geq S^{-}_{{\bar{P}}^{-1}}(\omega_1)$ by Lemma 2.5.
Then by (3.8)-(3.10), (4.18),  Lemma 2.5 and Lemma 2.7, we obtain
\bea i(y,\, 1) + 2S^+(y) - \nu(y,\,1)&=&i(\gamma_y)+2S^+_{\gamma_y(\tau)}(1)-\nu(\gamma_y)\nn\\
&=&i(\gamma_{\bar{P}}^{k})+2S^+_{{({\bar{P}}^{-1}M)}^{k}}(1)-\nu(\gamma_{\bar{P}}^{k})\nn\\&=&
\sum_{\omega^{k}=1} (i_\omega^{\bar{P}}(\gamma)+2S^+_{{\bar{P}}^{-1}M}(\omega)-\nu_{\omega}({\bar{P}}^{-1}M)),\nn\eea
which together with (4.25)-(4.26) and (4.20) implies
\bea i(y,\, 1) + 2S^+(y) - \nu(y,\,1)&\geq&
\sum_{\omega^{k}=1, \omega\neq 1} (\nu_\omega({\bar{P}}^{-1})+2S^+_{{\bar{P}}^{-1}M}(\omega)-\nu_{\omega}({\bar{P}}^{-1}M))
\nn\\&+&\sum_{i=1}^rS^{-}_{{\bar{P}}^{-1}M}(\alpha_{i})-
\sum_{i=1}^rS^{+}_{{\bar{P}}^{-1}M}(\alpha_i)+S^{-}_{{\bar{P}}^{-1}M}(\omega_1)-S^{+}_{{\bar{P}}^{-1}M}(1)\nn\\
&+&2S^+_{{\bar{P}}^{-1}M}(1)-\nu_{1}({\bar{P}}^{-1}M)\nn\\
&=&2n+\sum_{\omega^{k}=1, \omega\neq 1} (2S^+_{{\bar{P}}^{-1}M}(\omega)-\nu_{\omega}({\bar{P}}^{-1}M))+\sum_{i=1}^rS^{-}_{{\bar{P}}^{-1}M}(\alpha_{i})
\nn\\&-&
\sum_{i=1}^rS^{+}_{{\bar{P}}^{-1}M}(\alpha_i)+S^{-}_{{\bar{P}}^{-1}M}(\omega_1)
+S^+_{{\bar{P}}^{-1}M}(1)-\nu_{1}({\bar{P}}^{-1}M)\nn\\
&\geq&2n-\{\sum_{\omega^{k}=1, \omega\neq 1} (\nu_{\omega}({\bar{P}}^{-1}M)-2S^+_{{\bar{P}}^{-1}M}(\omega))
+\sum_{i=1}^rS^{+}_{{\bar{P}}^{-1}M}(\alpha_i)
\nn\\&+&(\nu_{1}({\bar{P}}^{-1}M)-S^+_{{\bar{P}}^{-1}M}(1))\}\nn\\&\geq&2n-n=n,\nn\eea
where we also used the fact that
\bea \sum_{\omega^{k}=1, \omega\neq 1} (\nu_{\omega}({\bar{P}}^{-1}M)-2S^+_{{\bar{P}}^{-1}M}(\omega))
+\sum_{i=1}^rS^{+}_{{\bar{P}}^{-1}M}(\alpha_i)
+(\nu_{1}({\bar{P}}^{-1}M)-S^+_{{\bar{P}}^{-1}M}(1))\leq n,\eea
in fact, by Lemma 2.5, we can get (4.27) if ${\bar{P}}^{-1}M$ is any basic form listed in (2.5)-(2.8)
and then (4.27) holds for any symplectic matrix ${\bar{P}}^{-1}M$ by Lemma 2.4 and Lemma 2.5. The proof is complete.\hfill\hb

{\bf Remark 4.4.} In Definition 1.1 of \cite{GM}, for the first time, V. Ginzburg and L. Macarini introduced the
important notion of {\it strong dynamical convexity} for contact forms invariant under a group action, supporting the standard contact structure on the
sphere, and in Theorem 1.6 of the same paper, they proved that any compact convex hypersurface which is symmetric with respect
to the origin satisfies this strong dynamical convexity condition. Comparing our Proposition 4.3 with Definition 1.1 of \cite{GM},
we can easily see that the compact convex P-cyclic symmetric hypersurfaces correspond to strong dynamical convex contact forms,
so our Proposition 4.3 gives a lot of examples satisfying this strong dynamical convexity condition, which also shows the validity
of the definition of {\it strong dynamical convexity}.

Now we can give the proofs of our main results:

{\bf Proof of Theorem 1.1.}

By Proposition 4.2, we denote the elements in $\T(\Sigma, \alpha)$ by
\bea \T(\Sigma, \alpha)=\{(\tau_j,y_j)\mid 1\leq j\leq s_1\}\cup \{(\tau_l,y_l), (\tau_l,Py_l)\mid s_1+1\leq l\leq s_1+s_2\},\nn\eea
where $\mathcal{O}(y_j)=\mathcal{O}(Py_j)$ for $1\leq j\leq s_1$, and $\mathcal{O}(y_l)\cap\mathcal{O}(Py_l)=\emptyset$ for $s_1+1\leq l\leq s_1+s_2$.
Since we have assumed $^\#\T(\Sigma, \alpha)=S$, then\bea S=s_1+2s_2.\eea
According to Lemma 3.6, we get an
injection map $p = p(\Sigma,\alpha) :\N\rightarrow \mathcal{T}(\Sg, \alpha)\times\N$. From (3.3), (4.13)-(4.14) and Lemma 3.6,
we can further require that\bea im(p)\subseteq \{(\tau_j,y_j)\mid 1\leq j\leq s_1+s_2\}\times \N.\nn\eea

By (3.7), we have $\hat i(y_j)>2$ for $1\le j\le s_1+s_2$, then we can use
the common index jump theorem (Theorems 4.3 and 4.4 of
\cite{LoZ1}, Theorems 11.2.1 and 11.2.2 of \cite{Lon4}) to obtain infinitely many
$(T, m_1,\ldots,m_{s_1+s_2})\in\N^{s_1+s_2+1}$ such that the following hold for every $j\in \{1,\ldots,s_1+s_2\}$:
\bea
\nu(y_j,\, 2m_j-1) &=&\nu(y_j,\, 1), \\
i(y_j,\, 2m_j) &\ge& 2T-\frac{e(\gamma_j(\tau_j))}{2}\ge 2T-n, \\
i(y_j,\, 2m_j)+\nu(y_j,\, 2m_j) &\le& 2T+\frac{e(\gamma_j(\tau_j))}{2}-1\le 2T+n-1, \\
i(y_j,\, 2m_j+1) &=& 2T+i(y_j,\,1), \\
i(y_j,\, 2m_j-1)+\nu(y_j,\, 2m_j-1)
 &=& 2T-(i(y_j,\,1)+2S^+(y_j)-\nu(y_j, 1)),
\eea
where $e(M)$ is the elliptic height of $M$ defined in \S1.
Note that (4.31) holds by Theorem 4.4 of \cite{LoZ1}, other parts
follows by Theorem 4.3 of \cite{LoZ1}.

By Lemma 2.4, we can assume $\gamma_j(\tau_j)$ can be connected within
$\Omega^0(\gamma_j(\tau_j))$ to
\be N_1(1,1)^{\diamond p_{j, -}}\diamond I_2^{\diamond p_{j, 0}}\diamond N_1(1,-1)^{\diamond  p_{j, +}}\diamond  G_j,
\qquad1\le j\le s_1+s_2,\ee
for some nonnegative integers $p_{j, -}$, $p_{j, 0}$, $p_{j, +}$ and some symplectic
matrix $G_j$ satisfying $1\not\in \sigma(G_j)$.
By (4.6), (4.34) and Lemma 2.5 we obtain
\be   2S^+(y_j)  = 2(p_{j, -} +p_{j, 0}) \ge 2,\quad1\le j\le s_1+s_2. \ee

By Corollary 3.2, we have
\be i(y_j,\,1)\ge n,\qquad 1\le j\le s_1+s_2.\ee
By (4.29), (4.33), (4.35) and (4.36) we have
\be i(y_j,\, 2m_j-1)  = 2T-(i(y_j,\,1)+2S^+(y_j))\le 2T-n-2.\ee
Combining (3.11)-(3.12) with (4.37),  for $m\ge 2$ we have
\bea i(y_j,{2m_j-m}) + \nu(y_j, {2m_j-m}) - 1
&\le& i(y_j,{2m_j-m+1}) - 1   \nn\\
&\le& i(y_j,{2m_j-1}) - 1   \nn\\
&\le& 2T - n - 3.  \eea

By (3.11), (4.32), and (4.36), for all $m\ge 2$ we obtain
\be  i(y_j,{2m_j+m}) > i(y_j,2m_j+1) = 2T + i(y_j,1)  \ge 2T+n. \ee

For every $1\le i\le n$, Denote by $p(T-i+1)=((\tau_{\rho(i)}, y_{\rho(i)}),\lambda(i))$, where $1\le i\le n$,
$\rho(i)\in\{1,\cdots, s_1+s_2\}$ and $\lambda(i)\in\N$. By the definition of $p$ and (3.13) we have
\bea i(y_{\rho(i)},\lambda(i))\leq 2T-2i+n\leq i(y_{\rho(i)},\lambda(i))+\nu(y_{\rho(i)},\lambda(i))-1,
\eea
which together with (4.38)-(4.39) yields
\be \lambda(i)\in\{2m_{\rho(i)}-1,\;2m_{\rho(i)}\}, \qquad \forall\;1\le i\le n.  \ee

{\bf Claim 1.} {\it If $y_{\rho(i)}$ is P-cyclic symmetric, then $\lambda(i)=2m_{\rho(i)}$. }

In fact, by (4.33) and (4.15) of Proposition 4.3, we have
\bea i(y_{\rho(i)},{2m_{\rho(i)}-1}) + \nu(y_{\rho(i)},{2m_{\rho(i)}-1}) - 1 \le 2T-n-1.\nn\eea
Thus Claim 1 holds by (4.40) and (4.41).

Since the map $p$ is injective, then by Claim 1 we have
\bea ^\#\{i\in \{1,\ldots,n\}\mid \rho(i)\leq s_1\}\leq s_1,\eea
and by (4.41) there holds
\bea ^\#\{i\in \{1,\ldots,n\}\mid \rho(i)> s_1\}\leq 2s_2.\eea
Therefore, combining (4.28) with (4.42)-(4.43) we obtain\bea ^\#\T(\Sigma)=^\#\T(\Sigma, \alpha)=S=s_1+2s_2\geq n.\nn\eea
The proof is complete.\hfill\hb

{\bf Proof of Theorem 1.2.}

Based on the proof of Theorem 1.1, we have the following claims:

{\bf Claim 2.} {\it If $\lambda(i)=2m_{\rho(i)}-1$, then $y_{\rho(i)}$ is not P-cyclic symmetric and
non-hyperbolic. }

The first statement follows directly from Claim 1. We prove the latter.

In fact, suppose $y_{\rho(i)}$ for some $i\in \{1, \ldots,n\}$ is hyperbolic. Then by (4.6), (4.36) and Lemma 2.4 we have $\nu(y_{\rho(i)})=S^+(y_{\rho(i)})=1$ and
\bea    i(y_{\rho(i)},\, 1)+2S^+(y_{\rho(i)})-\nu(y_{\rho(i)},\,1) = i(y_{\rho(i)},\, 1)+ 1 \ge n+1,\nn\eea
which together with (4.33) yields
\bea i(y_{\rho(i)},{2m_{\rho(i)}-1})+\nu(y_{\rho(i)},{2m_{\rho(i)}-1})-1 \le 2T-n-2 < 2T+n -2i, \nn\eea
which contradicts to (4.40). Thus Claim 2 holds.

{\bf Claim 3.} {\it When $n$ is even and $\lambda(i)=2m_{\rho(i)}$, $y_{\rho(i)}$ is non-hyperbolic. }

Suppose $y_{\rho(i)}$ is hyperbolic. Then we have $\nu(y_{\rho(i)})=1$ and
$e(\gamma_{\rho(i)}(\tau_{\rho(i)}))=2$. By (4.30), (4.31), correspondingly
we have
\be i(y_{\rho(i)},\, 2m_{\rho(i)})+\nu(y_{\rho(i)},\, 2m_{\rho(i)})-1\leq 2T-1\leq i(y_{\rho(i)},\, 2m_{\rho(i)}).\ee
On the other hand, by (4.40) we have  \bea i(y_{\rho(i)},2m_{\rho(i)})\leq 2T-2i+n\leq i(y_{\rho(i)},2m_{\rho(i)})+\nu(y_{\rho(i)},2m_{\rho(i)})-1,
\eea
which contradicts to (4.44), because $n$ is even. Hence
Claim 3 holds.

{\bf Claim 4.} {\it When $n$ is odd, let
$$  \I = \{i\in \{1, \ldots, n\}\;|\;\lambda(i)=2m_{\rho(i)}\;\;{\rm holds\;for}\;y_{\rho(i)}^{\lambda(i)}\}. $$
Then there exists at most one $i\in\I$ such that $y_{\rho(i)}$ is hyperbolic. Here we do not require
specially $y_{\rho(i)}$ is P-cyclic symmetric or not.}

In fact, suppose $y_{\rho(i)}$ is hyperbolic for some $i\in \I$, and then by (4.44)-(4.45)
we must have $2i = n+1$. Assume $y_{\rho(j)}$ is also hyperbolic for some $j\in \I\bs\{i\}$.
Then we obtain $2T-2j\neq 2T-2i = 2T-n-1$. Thus (4.44)-(4.45) with $i$ replaced by
$j$ imply $y_{\rho(j)}$ can not be hyperbolic. This completes the proof of Claim 4.

Now when $n$ is even, by Claim 2 and Claim 3, the closed characteristics found in Theorem 1.1 are all non-hyperbolic and thus
there exist at least $n$ non-hyperbolic closed characteristics.

When $n$ is odd, by Claim 2 and Claim 4, there exists at most one of the closed characteristics found in Theorem 1.1 is hyperbolic
and thus there exist at least $n-1$ non-hyperbolic closed characteristics. The proof is complete.\hfill\hb

{\bf Proof of Theorem 1.3.}

Based on the proof of Theorem 1.1,
without of loss of generality,  we denote the elements in $\{(\tau_l,y_l) \mid s_1+1\leq l\leq s_1+s_2\}$ by
\bea \{(\tau_l,y_l) \mid s_1+1\leq l\leq s_1+s_2\}&=&\{(\tau_j,y_j)\mid s_1+1\leq l\leq s_1+s_3\}\cup\nn\\&& \{(\tau_l,y_l), (\tau_l,P^2y_l)\mid s_1+s_3+1\leq l\leq s_1+s_3+s_4\},\eea
where $\mathcal{O}(y_l)=\mathcal{O}(P^2y_l)$ for $s_1+1\leq l\leq s_1+s_3$ and $\mathcal{O}(y_l)\cap\mathcal{O}(P^2y_l)=\emptyset$ for $s_1+s_3+1\leq l\leq s_1+s_3+s_4$.
Then by (4.46) we have\bea
s_2=s_3+2s_4.\eea

Since $\Sigma$ is also $P^2$-cyclic symmetric, then $\Phi(P^2u)=\Phi(u)$ for any
$u\in L_0^{\frac{\alpha}{\alpha-1}}(S^1,\R^{2n})$ from (4.13), and (4.14) with $P$ replaced by $P^2$ also holds, thus
by Lemma 3.6,
we can further require that\bea im(p)\subseteq \{(\tau_l,y_l)\mid 1\leq l\leq s_1+s_3+s_4\}\times \N.\nn\eea
Then by the same proof of Theorem 1.1, we obtain \bea ^\#\T(\Sigma)=^\#\T(\Sigma, \alpha)=s_1+2s_2\geq s_1+2(s_3+s_4)\geq n.\eea
On the other hand, by assumption we have \bea S=^\#\T(\Sigma)=^\#\T(\Sigma, \alpha)=n,\nn\eea
which together with (4.47)-(4.48) implies $s_4=0$, thus by (4.46)-(4.47) we have
\bea \mathcal{O}(y_j)=\mathcal{O}(P^2y_j), \quad s_1+1\leq j\leq s_1+s_2.\eea

Now we proceed our proof in two case according to the parity of $k$ which satisfies $P^k=I_{2n}$ and $k\geq3$.

{\bf Case 1.} $k$ is odd.

In this case, by (4.49) we have
\bea \mathcal{O}(y_j)=\mathcal{O}(P^2y_j)=\cdots=\mathcal{O}((P^2)^{\frac{k+1}{2}}y_j)=\mathcal{O}(Py_j), \quad s_1+1\leq j\leq s_1+s_2.\nn\eea
Because $\mathcal{O}(y_j)=\mathcal{O}(Py_j)$ for $1\leq j\leq s_1$, we get all the closed characteristics are P-cyclic symmetric.

{\bf Case 2.} $k$ is even.

In this case, by (4.49), for $s_1+1\leq j\leq s_1+s_2$, the  closed characteristics $(\tau_j, y_j)$ are $P^2$-cyclic symmetric.  Note that
we have $ker((P^2)^l-I_{2n})=0$ for $1\leq l<\frac{k}{2}$ since $ker(P^l-I_{2n})=0$ for $1\leq l<k$, and
$(P^2)^{\frac{k}{2}}=I_{2n}$ for $\frac{k}{2}\geq2$ since $k\geq3$ and $k$ is even, thus we can use
Proposition 4.3 with $P$ replaced by $P^2$ to obtain\bea    i(y_j,\, 1) + 2S^+(y_j) - \nu(y_j,\,1) \ge n, \quad s_1+1\leq j\leq s_1+s_2. \nn\eea
Then by the proof of Claim 1, we have $\lambda(i)=2m_{\rho(i)}$ for all $1\leq i\leq n$, which allow us to
use the same proof of Theorem 1.1 and obtain \bea ^\#\T(\Sigma)=^\#\T(\Sigma, \alpha)\geq s_1+s_2\geq n.\eea
On the other hand, by assumption we have \bea S=^\#\T(\Sigma)=^\#\T(\Sigma, \alpha)=n,\nn\eea
which together with (4.28) and (4.50) implies $s_2=0$, we get all the closed characteristics are P-cyclic symmetric. The proof is complete.\hfill\hb

\bibliographystyle{abbrv}

\end{document}